\documentclass[12pt]{article}
\usepackage{latexsym,amsfonts,amssymb}
\usepackage{indentfirst}

\usepackage[margin=3cm]{geometry}

\def\y{\begin{eqnarray*}}
\def\bd{\begin{description}}
\def\ey{\end{eqnarray*}}
\def\ebd{\end{description}}

\def\R{\mathbb{R}}

\begin{document}

\setlength{\baselineskip}{16pt}

\title{Periodic  solutions for neutral evolution equations with  delays
\thanks{Research supported by NNSF of China (11261053) and NSF of Gansu Province (1208RJZA129).
}}
\author{Qiang Li $^{1}$\thanks{
 \emph{Corresponding author.}~E-mail address: lznwnuliqiang@126.com (Q. Li), liyx@nwnu.edu.cn (Y. Li), xueye0528@126.com (H. Zhang).}, Yongxiang Li $^{2}$, Huanhuan Zhang $^{2}$
}

\date{\small
\begin{flushleft}
$1.$ Department of Mathematics, Shanxi Normal University, Linfen 041000, Peoples's Republic of China,
\vskip1mm
\noindent
$2.$ Department of Mathematics, Northwest Normal University, Lanzhou 730070, Peoples's Republic of China
\end{flushleft}
}
\maketitle

\begin{abstract}
\setlength{\baselineskip}{14pt}
 The aim is to study the  periodic solution problem for neutral evolution equation
  $$(u(t)-G(t,u(t-\xi)))'+Au(t)=F(t,u(t),u(t-\tau)),\ \ \ \ t\in\R$$in Banach space $X$, where $A:D(A)\subset X\rightarrow X$ is a closed linear operator, and $-A$ generates a compact analytic operator semigroup $T(t)(t\geq0)$.
  With the aid of the analytic operator semigroup theories and some fixed point theorems, we obtain the existence and uniqueness  of periodic mild solution for neutral evolution equations. The regularity of periodic mild solution for evolution equation with delay is  studied, and some the existence results of the classical and strong solutions are obtained.  In the end, we give an example to illustrate the applicability of abstract results. Our works greatly improve and generalize the relevant results of existing literatures.

\vspace{8pt}


\noindent {\bf  MR(2000) Subject Classification:\ } 34K30; 47H07; 47H08.

\end{abstract}

\section{Introduction}

The theory of partial differential equations with delays has extensive physical background and realistic mathematical model, and it has undergone a rapid development in the last fifty years see \cite{Hale1993,Wu1996} and references therein. More recently researchers have given special attentions to the study of equations in which the delay argument occurs in the derivative of the state variable as well as in the independent variable, so-called neutral differential equations. Neutral differential equations have many applications. It can model a lot of problems arising from engineering, such as population dynamics, transmission line, immune response or distribution of albumin in the blood etc.

In \cite{Wu1996,Wu1999} the authors studied a the partial neutral functional differential-difference equations which  is defined on a unit circle $S^{1}$:
\begin{eqnarray*}\frac{\partial}{\partial t}(u(x,t)-qu(x,t-\tau))&=&k\frac{\partial^{2}}{\partial x^{2}}(u(x,t)-qu(x,t-\tau))-au(x,t)-aqu(x,t-\tau)\\
&&-g(u(x,t)-qu(x,t-\tau)),\ \  x\in S^{1}, t\in \R,\qquad\qquad (1.1)\end{eqnarray*}
where $a,k,q$ are positive constants, $g:\R\rightarrow\R$ is continuously differentiable, $\tau\geq0$ which denotes the time delay.
Thereafter, more results on partial neutral functional differential equations are published, and we refer readers to \cite{Adimy2004,Adimy2006,Ezzinbi2004,Ezzinbi2010,Hernandeza2011,Regan2011}.The
idea of studying partial neutral functional differential equations with operators
satisfying Hille-Yosida condition, begins with \cite{Adimy1998}, where the authors studied the
following class of equation in a Banach space $X$.
$$(u(t)-Bu(t-\xi))'=A(u(t)-Bu(t-\xi))+F(u(t),u(t-\tau)),\eqno(1.2)$$
where $A$ satisfies the Hille-Yosida condition, $B$ are bounded linear operators
from $X$ into $X$, $F:X\times X\rightarrow X$ is continuous, $\xi,\tau$ are poeitive constants which denote the time delays. It has
been proved in particular, that the solutions generate a locally Lipschitz continuous integrated semigroup.

The problems concerning periodic solutions of partial neutral functional differential equations with delay are an important area of investigation since they can take into account seasonal fluctuations occurring in the phenomena appearing in the models, and have been studied by some researchers in recent years. Specially, the existence of periodic solutions of neutral evolution equations with delay has been considered by several authors, see  \cite{Babram1996,Hernandez1998,Benkhalti2006,Fu2007,Fu2008,Benkhalti2010,Ezzinbi2014}. For the delayed evolution equations without neutral term,
 the existence of periodic solutions has been discussed by more
 authors, see \cite{Burton1991,Xiang1992,Liu1998,Liu2000,Liu2003} and references therein.
Naturally,
fixed point theorems play a significant role in the investigation of the existence of
periodic solutions. It is well known that  the
Massera's approach (see \cite{Massera1950}) on periodic partial functional differential equations explains the relationship between the existence of bounded solutions and the existence of periodic solutions. However, in many of the studies mentioned above, the key assumption of prior
boundedness of solutions was employed and the most important feature is to show that Poincar\'{e}'s mapping
$$P_{\omega}(\phi)=u_{\omega}(\cdot,\phi),$$
is condensing, where $\omega$ is a period of the system and $u$ the unique mild solution determined
by $\phi$. Therefore,  a fixed point theorem can be used to derive periodic solutions.

Recently, Zhu, Liu and Li in \cite{Zhu2008} investigated the existence of time periodic solutions for a class of
one-dimensional parabolic evolution equation with delays. They obtained the existence of time periodic solutions by constructing some suitable Lyapunov functionals and establishing the
prior bound for all possible periodic solutions. And, Li in \cite{Li2011} discussed the existence of the time periodic solution for the evolution
equation with multiple delays in a Hilbert space $H$
$$u'(t)+Au(t)=F(t,u(t),u(t-\tau_{1},\cdots,u(t-\tau_{n}))),\ \ \ \ t\in\R,\eqno(1.3)$$
where $A : D(A)\subset H\rightarrow H$ is a positive definite selfadjoint operator, $F : \R \times H^{n+1} \rightarrow H$ is a
nonlinear mapping which is $\omega$-periodic in $t$, and $\tau_{1},\tau_{2},\cdots,\tau_{n}$ are positive constants which
denote the time delays. By using periodic extension and Schauder fixed point theorem, the author presented essential conditions  on the nonlinearity $F$ to guarantee that the equation has $\omega$-periodic solutions.

Motivated by the papers mentioned above, the aim of this work is to study the existence of periodic solution for some the partial neutral functional differential equations. Our discussion will be made in a frame of abstract Banach spaces.

Throughout this paper, $X$ is a Banach space provided with norm $\|\cdot\|$ and $A:D(A)\subset X\rightarrow X$ is a closed linear operator, and $-A$ generates a compact analytic operator semigroup $T(t)(t\geq0)$ in Banach space $X$.
Let $G,F $ be appropriate continuous
functions which will be specified later, and  $G(t,\cdot),F(t,\cdot,\cdot)$ be $\omega$-periodic in $t$.

Under the above assumptions  we discuss the existence and uniqueness of $\omega$-periodic solutions of the
abstract neutral functional differential equations with delays in $X$ of the form
$$(u(t)-G(t,u(t-\xi)))'+Au(t)=F(t,u(t),u(t-\tau)),\ \ \ \ t\in\R,\eqno(1.4)$$
where $\xi,\tau$ are positive constants which denote the time delays. The purpose of the present note is to extend and develop the work in \cite{Zhu2008,Li2011}, that is, we will discuss the existence  and regularity of periodic solutions for Eq. (1.4). The obtained results will also improve the main results in \cite{Hernandez1998,Benkhalti2010,Ezzinbi2014}. In this paper, it is worth mentioning that assumption of prior boundedness of solutions is not necessary.  More precisely, the nonlinear term $F$ only satisfies some  growth
conditions and the function $G$ and $F$ may not be defined on the whole space $X$. These conditions are much weaker than Lipschitz conditions.

The paper is organized as follows. In Section 2, we  collect some known notions and results on the fractional powers of the generator of an analytic semigroup and  provide
preliminary results to be used in theorems stated and proved in the paper.
In Section 3, we apply the operator semigroup theory to find the periodic mild solutions for  Eq.(1.4) and in Section 4, we investigate conditions for Eq.(1.4) to have the calssical and strong periodic solutions. In the last section, we give an example to illustrate the applicability of abstract results obtained in Section 3 and Section 4.

\section{Preliminaries}

Throughout this paper, we assume that $X$ is a Banach space with norm $\|\cdot\|$,
 that $A:D(A)\subset X\rightarrow X$ is a closed linear operator and  $-A$ generates a compact analytic operator semigroup $T(t)(t\geq0)$ in Banach space $X$. For the theory of semigroups of linear operators we refer to \cite{Pazy1993}.

 We only recall here some notions and properties that are essential for us.
For a general $C_{0}$-semigroup $T(t)(t\geq0)$, there exist $M\geq1$ and $\nu\in \R$ such that (see \cite{Pazy1993})
$$\|T(t)\|\leq Me^{\nu t},\quad t\geq0.\eqno(2.1)$$
Let
$$\nu_{0}=\inf\{\gamma\in \R |\ \mathrm{There\ exists}\ M\geq1\ \mathrm{ such \ that}\  \|T(t)\|\leq Me^{\nu t},\ \forall t\geq0\},$$
 then $\nu_{0}$ is called the growth exponent of the semigroup $T(t)(t\geq0)$. Furthermore, $\nu_{0}$ can be also obtained by the following formula
$$\nu_{0}=\limsup\limits_{t\rightarrow +\infty}\frac{\ln\|T(t)\|}{t}.\eqno(2.2)$$
If $C_{0}$-semigroup $T(t)$ is analytic on $(0,+\infty)$, it is well known that $\nu_{0}$ can also be determined by $\sigma(A)$ (see\cite{Pazy1993,Triggiani1975}),
$$\nu_{0}=-\inf\{\mathrm{Re} \lambda |\ \lambda\in\sigma(A)\},\eqno(2.3)$$
where $-A$ is the infinitesimal generator of $C_{0}$-semigroup $T(t)(t\geq0)$. We know that $T(t)(t\geq0)$ is continuous in the uniform operator topology for $t>0$ if $T(t)(t\geq0)$ is compact semigroup or analytic semigroup(see \cite{ Hino1987}).

In particular, if $T(t)(t\geq0)$ is analytic semigroup with infinitesimal generator $A$ satisfying $0 \in \rho(A)$($\rho(A)$ is the resolvent set of $A$).
Then for any $\alpha > 0$, we can define $A^{-\alpha}$ by
$$A^{-\alpha}:=\frac{1}{\Gamma(\alpha)}\int_{0}^{\infty}t^{\alpha-1}T(t)dt.$$
It follows that each $A^{- \alpha}$ is an injective continuous endomorphism of $X$. Hence we
can define $A^{\alpha}$ by $A^{\alpha} := (A^{-\alpha})^{-1}$,
 which is a closed bijective linear operator in $X$. Furthermore, the subspace $D(A^{\alpha})$ is
 dense in $X$ and the expression
$$\|x\|_{\alpha}:=\|A^{\alpha}x\|,\ \ \ x\in D(A),$$
 defines a norm on $D(A^{\alpha})$.  Hereafter we
 respresent $X_{\alpha}$ to the space $D(A^{\alpha})$ endowed with the norm $\|\cdot\|_{\alpha}$ and denote
  by $C_{\alpha}$ the operator norm of $A^{-\alpha}$, i.e., $C_{\alpha}:=\|A^{-\alpha}\|$. The following
  preperties are well known(\cite{Pazy1993}).
 \vskip3mm
\noindent\textbf{Lemma 2.1.} \emph{ If $T(t)(t\geq0)$ is
 analytic semigroup with infinitesimal generator $A$ satisfying $0 \in \rho(A)$, then\\
(i) $D(A^{\alpha})$ is a Banach space for $0\leq\alpha\leq1$;\\
(ii) $A^{-\alpha}$ is a bounded linear operator for $0\leq\alpha\leq1$ in $X$;\\
(iii) $T(t):X\rightarrow D(A^{\alpha})$ for each $t>0$;\\
(iv) $A^{\alpha}T(t)x=T(t)A^{\alpha}x$ for each $x\in D(A^{\alpha})$ and $t\geq0$;\\
(v) for every $t > 0$, $A^{\alpha}T(t)$ is bounded in $X$ and there exists $M_{\alpha} > 0$ such that
$$\|A^{\alpha}T(t)\|\leq M_{\alpha}t^{-\alpha},$$
moreover, if $\alpha\in(0,1)$ then $M_{\alpha}=M\Gamma(\alpha)$, where $M$ is given by (2.1);\\
(vi)  $X_{\beta}\hookrightarrow X_{\alpha}$ for $0 \leq \alpha \leq\beta\leq1$ (with $X_{0}=X$ and $X_{1}=D(A)$), and the embedding is continuous.
and the embedding $X_{\beta} \hookrightarrow X_{\alpha}$
is compact whenever the resolvent operator of $A$ is compact.}
\vskip2mm

Observe by Lemma 2.1 (iii) and (iv) that the restriction
 $T_{\alpha}(t)$ of $T(t)$ to $X_{\alpha}$ is exactly the
  part of $T(t)$ in $X_{\alpha}$. Moreover, for any $x\in X_{\alpha}$, we have
$$\|T_{\alpha}(t)x\|_{\alpha}=\|A^{\alpha}T(t)x\|=\|T(t)A^{\alpha}x\|\leq\|T(t)\|\cdot\|A^{\alpha}x\|=\|T(t)\|\cdot\|x\|_{\alpha},\eqno(2.4)$$
and
$$\|T_{\alpha}(t)x-x\|_{\alpha}=\|A^{\alpha}T(t)x-A^{\alpha}x\|=\|T(t)A^{\alpha}x-A^{\alpha}x\|\rightarrow 0,\quad t\rightarrow 0,\eqno(2.5)$$
it follows that $T_{\alpha}(t)(t\geq0)$ is a strongly continuous semigroup on $X_{\alpha}$ and $\|T_{\alpha}(t)\|_{\alpha}\leq \|T(t)\|$ for all $t \geq 0$. To prove our main results, we need the following lemmas.
\vskip3mm
\noindent\textbf{Lemma 2.2.} (\cite{Liu2009}) \emph{If $T(t)(t\geq0)$ is a compact semigroup in $X$, then $T_{\alpha}(t)(t\geq 0)$ is a compact semigroup in $X_{\alpha}$.}
\vskip3mm
\noindent\textbf{Lemma 2.3.} (\cite{Chang2009}) \emph{If $X$ is reflexive, then $X_{\alpha}$ is also reflexive.}
\vskip3mm

Now, recall some basic facts on abstract linear evolutions, which are needed to prove our main results.

Let $J$ denote the infinite interval $[0,\infty)$ and   $h:J\rightarrow X$, consider the initial value problem of the linear evolution equation
$$\left\{\begin{array}{ll}
u'(t)+Au(t)=h(t),\ t\in J,\\[8pt]
u(0)=x_{0}.
 \end{array} \right.\eqno(2.6)$$
 It is well known, when $x_{0}\in D(A)$ and $h\in C^{1}(J,X)$,
 the initial value problem (2.6) has a unique classical solution
 $u\in C^{1}(J,X)\cap C(J,X_{1})$ expressed by
 $$u(t)=T(t)x_{0}+\int^{t}_{0}T(t-s)h(s)ds. \eqno(2.7)$$
Generally, for $x_{0}$ and $h\in C(J,X)$,
the function $u$ given by (2.7) belongs to $C(J,X)$ and it is
called a mild solution of the linear evolution equation (2.6). A mild solution $u$
of Eq. (2.6) is called a strong solution if $u$ is continuously differentiable a.e. on $J$, $u'\in L_{loc}^{1}(J,X)$ and satisfies Eq. (2.6). Furthermore, we have the following results.

\vskip3mm
\noindent\textbf{Lemma 2.4.}(\cite{Pazy1993}) \emph{ Let $h\in C([0,a],X)(a>0)$, $0\leq\alpha<\beta\leq1$, $\mu=\beta-\alpha$, $x_{0}\in X_{\beta}$, then the mild solution $u$ of Eq. (2.6) satisfies $u\in c^{\mu}([0,a],X_{\alpha})$.}

\vskip3mm
\noindent\textbf{Lemma 2.5.}(\cite{Pazy1993})\emph{ Let $h\in C^{\mu}([0,a],X)(a>0)$, $0<\mu<1$, $x_{0}\in X$, then the mild solution $u$ of Eq. (2.6) is classical solution on $[0,a]$ and satisfies
$$u\in C^{1}((0,a],X)\cap C([0,a],X_{1}).$$}
\indent Let $C_{\omega}(\R,X)$ denote the Banach space $\{u\in C(\R,X)|\ u(t)=u(t+\omega),t\in\R\}$ endowed the
maximum norm $\|u\|_{C}=\max_{t\in J_{0}}\|u(t)\|$,  and $C_{\omega}(\R,X_{\alpha})$ denote the Banach space $\{u\in C(\R,X_{\alpha})|\ u(t)=u(t+\omega),t\in\R\}$ endowed the
maximum norm $\|u\|_{C\alpha}=\max_{t\in J_{0}}\|u(t)\|_{\alpha}$. Clearly, $C_{\omega}(\R,X_{\alpha})\hookrightarrow C_{\omega}(\R,X)$.

Given $h\in C_{\omega}(\R,X)$, we consider the existence of $\omega$-periodic mild solution of linear evolution equation
$$u'(t)+Au(t)=h(t),\quad t\in \R.\eqno(2.8)$$
\noindent\textbf{Lemma 2.6.}(\cite{Li2005}) \emph{ If $-A$ generates an exponentially stable $C_{0}$-semigroup $T(t)(t\geq0)$ in $X$ that is $\nu_{0}<0$,
then for $h\in C_{\omega}(\R,X)$,
the linear evolution equation (2.8) exists a unique $\omega$-periodic mild
solution $u$, which can be expressed by
  $$u(t)=(I-T(\omega))^{-1}\int_{t-\omega}^{t}T(t-s)h(s)ds:=(Ph)(t),\eqno(2.9)$$ and the solution operator $P:C_{\omega}(\R,X)\rightarrow C_{\omega}(\R,X)$ is a bounded linear operator.}
\vskip3mm
\noindent\textbf{Proof.} For any $\nu\in (0,|\nu_{0}|)$, there exists $M>0$ such that
$$\|T(t)\|\leq Me^{-\nu t}\leq M, \quad t\geq0.$$
In $X$, define the equivalent norm $|\cdot|$ by
$$|x|=\sup\limits_{t\geq0}\|e^{\nu t}T(t)x\|,$$
then $\|x\|\leq |x|\leq M\|x\|$. By $|T(t)|$ we denote the norm of $T(t)$
in $(X,|\cdot|$), then for $t\geq0$, it is easy to obtain that
$|T(t)|<e^{-\nu t}$.
Hence, $(I-T(\omega))$ has bounded inverse operator
$$(I-T(\omega))^{-1}=\sum_{n=0}^{\infty}T(n\omega),$$
and its norm satisfies
$$|(I-T(\omega))^{-1}|\leq \frac{1}{1-|T(\omega)|}\leq \frac{1}{1-e^{-\nu \omega}}.\eqno(2.10)$$
Set $$x_{0}=(I-T(\omega))^{-1}\int^{\omega}_{0}T(t-s)h(s)ds:=Bh,\eqno(2.11)$$
 then the mild solution $u(t)$  of  the linear initial value problem  (2.6)
given by (2.7) satisfies the periodic boundary condition $u(0)=u(\omega)=x_{0}$.
For $t\in \mathbb{R}^{+}$, by (2.7) and the properties of the semigroup $T(t)(t\geq0)$,
we have
\begin{eqnarray*}
u(t+\omega)&=&T(t+\omega)u(0)+\int_{0}^{t+\omega}T(t+\omega-s)h(s)ds\\
&=&T(t)\Big(T(\omega)u(0)+\int_{0}^{\omega}T(\omega-s)h(s)ds\Big)+\int^{t}_{0}T(t-s)h(s-\omega)ds\\
&=&T(t)u(0)+\int^{t}_{0}T(t-s)h(s)ds\ =u(t).
\end{eqnarray*}
Therefore, the $\omega$-periodic extension of $u$ on $\R$ is a unique $\omega$-periodic
mild solution of Eq.(2.8). By (2.7) and (2.11), the $\omega$-periodic mild solution
can be  expressed by
\begin{eqnarray*}\qquad\qquad\ \
u(t)&=&T(t)B(h)+\int^{t}_{0}T(t-s)h(s)ds\\
&=&(I-T(\omega))^{-1}\int^{t}_{t-\omega}T(t-s)h(s)ds:=(Ph)(t).\qquad\qquad\ \qquad (2.12)
\end{eqnarray*}
Evidently, $P:C_{\omega}(\R,X)\rightarrow C_{\omega}(\R,X)$
is a bounded linear operator. In fact, for every $h\in C_{\omega}(\R,X)$
\begin{eqnarray*}\|Qh(t)\|&=&\Big\|(I-T(\omega))^{-1}\int^{t}_{t-\omega}T(t-s)h(s)ds\Big\|\\
&\leq&\|(I-T(\omega))^{-1}\|\cdot\int^{t}_{t-\omega}\|T(t-s)\|ds\|h\|_{C}\\[4pt]
&\leq& CM\omega\|h\|_{C},\end{eqnarray*}
where $C:=\|(I-T(\omega))^{-1}\|$, which implies that $P$ is bounded.  This completes the proof of Lemma 2.6. $\Box$

To prove our main results, we also need the following lemma.

\vskip3mm
\noindent\textbf{Lemma 2.7.}(\cite{Sadovskii1967}){\emph{ Assume that $Q$ is a condensing operator on a Banach space $X$, i.e.,$Q$ is continuous and takes bounded sets
into bounded sets, and $\alpha(Q(D))<\alpha(D)$  for every bounded set $D$ of $X$ with $\alpha(D)>0$. If $Q(\overline{\Omega})\subset\overline{\Omega}$ for a convex, closed, and bounded set
 $\overline{\Omega}$ of $X$, then $Q$ has a fixed point in $\overline{\Omega}$ (where $\alpha(\cdot)$ denotes the kuratowski measure of non-compactness).}}
 \vskip3mm
\noindent\textbf{Remark 2.8.} \emph{It is easy to see that, if $Q=Q_{1}+Q_{2}$ with $Q_{1}$ a completely continuous operator and $Q_{2}$ a contractive one, then $Q$ is a condensing operator on $X$.}

\section{Existence of Mild solution}
\vskip3mm
Now, we are in a position to state and prove our main results of this section.
\vskip1mm
\noindent\textbf{Theorem 3.1.}\emph{ Let $A:D(A)\subset X\rightarrow X$ be a closed linear operator, and $-A$ generate a compact and exponentially stable  analytic operator semigroup $T(t)(t\geq0)$ in Banach space $X$. For $\alpha\in [0,1)$, we assume that  $G:\R\times X_{\alpha}\rightarrow X_{1}$ and $F : \R \times X_{\alpha}^{2} \rightarrow X$ are continuous
functions, and for every $x,x_{0},x_{1}\in X_{\alpha}$, $G(t,x),F(t,x_{0},x_{1})$ are $\omega$-periodic in $t$. If the following conditions
\vskip1mm
\noindent (H1) for any $r>0$, there exists a positive value function $h_{r}:\R\rightarrow \R^{+}$
 such that
$$\sup_{\|x_{0}\|_{\alpha},\|x_{1}\|_{\alpha}<r}\|F(t,x_{0},x_{1})\|\leq h_{r}(t),\ \ \ t\in \R,$$
function $s\mapsto \frac{h_{r}(s)}{(t-s)^{\alpha}}$ belongs to $L_{loc}(\R,\R^{+})$  and there is a positive constant $\gamma>0$ such that
$$\liminf_{r\rightarrow\infty}\frac{1}{r}\int_{t-\omega}^{t}\frac{h_{r}(s)}{(t-s)^{\alpha}}ds=\gamma<\infty,\ \ \ t\in\R;$$
\vskip1mm
\noindent(H2) $G(t,\theta)=\theta$ for $t\in \R$, and there is a  constant  $L\geq0$ such that
$$\|AG(t,x)-AG(t,y)\|\leq L\|x-y\|_{\alpha},\ \ \ t\in \R, x,y\in X_{\alpha};$$
\vskip1mm
\noindent (H3) $CM_{\alpha}\gamma+C_{1-\alpha}L+CM_{\alpha}L\frac{\omega^{1-\alpha}}{1-\alpha}<1$, where $C=\|(I-T(\omega))^{-1}\|$,\vskip2mm
\noindent hold, then Eq. (1.4) has at least one $\omega$-periodic mild solution $u$.}
\vskip3mm

\noindent\textbf{Proof } From the assumption, we know that $G(t,u(t-\xi))\in D(A)$ for every $u\in C_{\omega}(\R,X_{\alpha})$, thus, we can  rewrite Eq.(1.4) as following
\begin{eqnarray*}\ \ \qquad\qquad&&(u(t)-G(t,u(t-\xi)))'+A(u(t)-G(t,u(t-\xi)))\\[8pt]
&=&F(t,u(t),u(t-\tau))-AG(t,u(t-\xi)), \ \ \ \quad t\in\R.\ \ \quad\qquad\qquad\qquad(3.1)\end{eqnarray*}

For any  $r>0$, let
$$\overline{\Omega}_{r}=\{u\in C_{\omega}(\R,X_{\alpha})\ |\  \|u\|_{C\alpha}\leq r\}.\eqno(3.2)$$
Note that $\overline{\Omega}_{r}$ is a closed ball in $C_{\omega}(\R,X_{\alpha})$ with centre $\theta$ and radius $r$. Moreover, by the condition (H2), it follows that,
\begin{eqnarray*}\|T(t-s)AG(s,u(s-\xi))\|_{\alpha}&=&\|A^{\alpha}T(t-s)A(G(s,u(s-\xi))-G(s,\theta))\|\\[8pt]
&\leq& \|A^{\alpha}T(t-s)\|\cdot \|A(G(s,u(s-\xi))-G(s,\theta))\|\\[8pt]
&\leq& \frac{M_{\alpha}L}{(t-s)^{\alpha}}\|u(s-\xi)\|_{\alpha},\end{eqnarray*}
which implies that $s\rightarrow T(t-s)AG(s,u(s-\xi))$ is integrable on $[t-\omega,t]$ for each $u\in \overline{\Omega}_{r}$.

Hence, we can define the operator $Q$ on $C_{\omega}(\R,X_{\alpha})$ by
\begin{eqnarray*}\qquad Qu(t)&:=&(I-T(\omega))^{-1}\int_{t-\omega}^{t}T(t-s)F(s,u(s),u(s-\tau))ds+G(t,u(t-\xi))\\
&&-(I-T(\omega))^{-1}\int_{t-\omega}^{t}T(t-s)AG(s,u(s-\xi))ds. \ \ t\in \R.\qquad\qquad(3.3)\end{eqnarray*}
From Lemma 2.6, it is sufficient to prove that $Q$ has a fixed point.

 Now, we  show that there is a positive constant $r$ such that $Q(\overline{\Omega}_{r})\subset\overline{\Omega}_{r}$. If this were not case, then for any $r>0$, there exist $u_{r}\in\overline{\Omega}_{r}$ and $t_{r}\in \R$ such that $\|Qu_{r}(t_{r})\|_{\alpha}>r$. Thus, we see by (H1),(H2) and  (H3) that
\begin{eqnarray*}
r&<&\|Qu_{r}(t_{r})\|_{\alpha}\\[8pt]
&\leq&\Big\|(I-T(\omega))^{-1}\int_{t_{r}
-\omega}^{t_{r}}T(t_{r}-s)F(s,u_{r}(s),u_{r}(s-\tau))ds\Big\|_{\alpha}
+\|G(t_{r},u_{r}(t_{r}-\xi))\|_{\alpha}\\
&&+\Big\|(I-T(\omega))^{-1}\int_{t_{r}-\omega}^{t_{r}}T(t_{r}-s)AG(s,u_{r}(s-\xi))ds\Big\|_{\alpha}\\
&\leq& \|(I-T(\omega))^{-1}\|\cdot\int_{t_{r}-\omega}^{t_{r}}\|A^{\alpha}T(t_{r}-s)\|\cdot \|F(s,u_{r}(s),u_{r}(s-\tau))\|ds\\[8pt]
&&+\|A^{\alpha-1}(AG(t_{r},u_{r}(t_{r}-\xi))-AG(t_{r},\theta))\|\\[8pt]
&&+\|(I-T(\omega))^{-1}\|\cdot\int_{t_{r}
-\omega}^{t_{r}}\|A^{\alpha}T(t_{r}-s)\|\cdot\|AG(s,u_{r}(s-\xi))-AG(s,\theta)\|ds\\[8pt]
&\leq& CM_{\alpha}\int_{t_{r}-\omega}^{t_{r}}\frac{h_{r}(s)}{(t_{r}-s)^{\alpha}}ds+C_{1-\alpha}L\|u_{r}\|_{C\alpha}+
CM_{\alpha}L\int_{t_{r}-\omega}^{t_{r}}\frac{1}{(t_{r}-s)^{\alpha}}ds\|u_{r}\|_{C\alpha}\\[8pt]
&\leq& CM_{\alpha}\int_{t_{r}-\omega}^{t_{r}}\frac{h_{r}(s)}{(t_{r}-s)^{\alpha}}ds
+C_{1-\alpha}Lr+CM_{\alpha}L\frac{\omega^{1-\alpha}}{1-\alpha}r. \end{eqnarray*}
Dividing on both sides by $r$ and taking the lower limit as $r\rightarrow\infty$, we have
$$CM_{\alpha}L\gamma+C_{1-\alpha}L+CM_{\alpha}L\frac{\omega^{1-\alpha}}{1-\alpha}\geq1,\eqno(3.4)$$
which contradicts (H3). Hence, there is a positive constant $r$ such that $Q(\overline{\Omega}_{r})\subset\overline{\Omega}_{r}$.

To show that the operator $Q$ has a fixed point on $\overline{\Omega}_{r}$,  we also
introduce the decomposition $Q=Q_{1}+Q_{2}$, where
$$Q_{1}u(t):=(I-T(\omega))^{-1}\int_{t-\omega}^{t}T(t-s)F(s,u(s),u(s-\tau))ds,\ \eqno(3.5)$$
$$Q_{2}u(t):=G(t,u(t-\xi))-(I-T(\omega))^{-1}\int_{t-\omega}^{t}T(t-s)AG(s,u(s-\xi))ds.\eqno(3.6)$$
Then we will prove that $Q_{1}$ is a compact operator and $Q_{2}$ is a contraction.

Firstly, we prove that $Q_{1}$ is a compact operator. Let $\{u_{n}\}\subset\overline{\Omega}_{r}$ with $u_{n}\rightarrow u$ in $\overline{\Omega}_{r}$, then by the continuity of $F$, we have
$$F(t,u_{n}(t),u_{n}(t-\tau))\rightarrow F(t,u(t),u(t-\tau)), \ \ n\rightarrow\infty,$$
for each $t\in\R$. Since $\|F(t,u_{n}(t),u_{n}(t-\tau))- F(t,u(t),u(t-\tau))\|\leq 2h_{r}(t)$ for all $t\in \R$, then the dominated convergence theorem ensure that
\begin{eqnarray*}
&&\|Q_{1}u_{n}(t)-Q_{1}u(t)\|_{\alpha}\\[8pt]
&=&\Big\|(I-T(\omega))^{-1}\int_{t-\omega}^{t}T(t-s)F(s,u_{n}(s),u_{n}(s-\tau))
-F(s,u(s),u(s-\tau))ds\Big\|_{\alpha}\\
&\leq&\|(I-T(\omega))^{-1}\| \int_{t-\omega}^{t}\|A^{\alpha}T(t-s)\|\cdot\|F(s,u_{n}(s),u_{n}(s-\tau))-F(s,u(s),u(s-\tau))\|ds\\
&\leq&CM_{\alpha}\int_{t-\omega}^{t}\frac{\|F(s,u_{n}(s),u_{n}(s-\tau))
-F(s,u(s),u(s-\tau))\|}{(t-s)^{\alpha}}ds\\[8pt]
&\rightarrow& 0 \ \ \mathrm{as}\ \  n\rightarrow\infty
\end{eqnarray*}
which implies that $\|Q_{1}u_{n}-Q_{1}u\|_{C\alpha}\rightarrow 0$ as $n\rightarrow\infty$, i.e. $Q_{1}$ is continuous.

It is easy to see  that $Q_{1}$ maps
 $\overline{\Omega}_{r}$ into a bounded set in $C_{\omega}(\R,X_{\alpha})$. Now, we demonstrate that $ Q(\overline{\Omega}_{r})$ is equicontinuous. For every
$u\in\overline{\Omega}_{r}$, by the periodicity of $u$, we only consider it on $[0,\omega]$.
Set $ 0\leq t_{1}<t_{2}\leq\omega$, we get that
\begin{eqnarray*}
&&Q_{1}u(t_{2})-Q_{1}u(t_{1})\\[8pt]
&=&(I-T(\omega))^{-1}\int^{t_{2}}_{t_{2}-\omega}T(t_{2}-s)F(s,u(s),u(s-\tau))ds\\
&~&-(I-T(\omega))^{-1}\int^{t_{1}}_{t_{1}-\omega}T(t_{1}-s)F(s,u(s),u(s-\tau))ds\\
&=&(I-T(\omega))^{-1}\int^{t_{1}}_{t_{2}-\omega}(T(t_{2}-s)-T(t_{1}-s))F(s,u(s),u(s-\tau))ds\\
&~&-(I-T(\omega))^{-1}\int^{t_{2}-\omega}_{t_{1}-\omega}T(t_{1}-s)F(s,u(s),u(s-\tau))ds\\
&~&+(I-T(\omega))^{-1}\int^{t_{2}}_{t_{1}}T(t_{2}-s)F(s,u(s),u(s-\tau))ds\\[8pt]
&:=&I_{1}+I_{2}+I_{3}.
\end{eqnarray*}
It is clear that
$$\|Q_{1}u(t_{2})-Q_{1}u(t_{1})\|_{\alpha}\leq \|I_{1}\|_{\alpha}+\|I_{2}\|_{\alpha}+\|I_{3}\|_{\alpha}.\eqno(3.7)$$
Thus, we only need to check $\|I_{i}\|_{\alpha}$ tend to $0$  independently of $u\in\overline{\Omega}_{r}$
when $t_{2}-t_{1}\rightarrow 0,i=1,2,3$. From  the continuity of $t \mapsto \|T(t)\|$ for $t > 0$ and the condition (H1), we can easily see
\begin{eqnarray*}
\|I_{1}\|_{\alpha}&\leq&C\cdot\int^{t_{1}}_{t_{2}-\omega}\|A^{\alpha}
(T(t_{2}-s)-T(t_{1}-s))\|\cdot\|F(s,u(s),u(s-\tau))\|ds\\
&\leq &  C\cdot\int^{t_{1}}_{t_{2}-\omega}\Big\|T\Big(\frac{t_{2}-s}{2}+\frac{t_{2}-t_{1}}{2}\Big)
-T\Big(\frac{t_{1}-s}{2}\Big)\Big\|\cdot\Big\|A^{\alpha}T\Big(\frac{t_{1}-s}{2}\Big)\Big\| \cdot h_{r}(s) ds\\[8pt]
&\leq&CM_{\alpha}\int^{t_{1}}_{t_{2}-\omega} \Big\|T\Big(\frac{t_{2}-s}{2}+\frac{t_{2}-t_{1}}{2}\Big)-T\Big(\frac{t_{1}-s}{2}\Big)\Big\|
\cdot\frac{h_{r}(s)}{(\frac{t_{1}-s}{2})^{\alpha}}ds\\[8pt]
&\rightarrow&0, \ \mathrm{as} \ t_{2}-t_{1}\rightarrow 0,\\[10pt]
\|I_{2}\|_{\alpha}&\leq&C\cdot\int^{t_{2}-\omega}_{t_{1}-\omega}\|AT(t_{1}-s)\|\cdot\|F(s,u(s),u(s-\tau))\|ds\\
&\leq &CM_{\alpha}\cdot\int^{t_{2}-\omega}_{t_{1}-\omega} \frac{h_{r}(s)}{(t_{1}-s)^{\alpha}}ds\\[8pt]
&\rightarrow&0,  \ \mathrm{as} \ t_{2}-t_{1}\rightarrow 0,\\[10pt]
\|I_{3}\|_{\alpha}&\leq&C\cdot\int^{t_{2}}_{t_{1}}\|A(T(t_{2}-s))\|\cdot\|F(s,u(s),u(s-\tau)))\|ds\\
&\leq &CM_{\alpha}\cdot\int^{t_{2}}_{t_{1}}\frac{h_{r}(s)}{(t_{2}-s)^{\alpha}}ds \\[8pt]
&\rightarrow&0,  \ \mathrm{as} \ t_{2}-t_{1}\rightarrow 0.
\end{eqnarray*}
As a result, $\| Q_{1}u(t_{2}) -Q_{1}u(t_{1})\|_{\alpha}$ tends to 0 independently of $u\in \overline{\Omega}_{r}$
as $t_{2}- t_{1}\rightarrow0$, which means that $Q_{1}(\overline{\Omega}_{r})$ is equicontinuous.

It remains to show that $(Q_{1}\overline{\Omega}_{r})(t)$ is relatively compact in $X_{\alpha}$ for all $t\in \R$.
To do this, we define a set $(Q_{\varepsilon}\overline{\Omega}_{r})(t)$ by
$$(Q_{\varepsilon}\overline{\Omega}_{r})(t):=\{(Q_{\varepsilon}u)(t)|\ u\in \overline{\Omega}_{r},\ 0<\varepsilon<\omega,\ t\in \R\}, \eqno(3.8)$$
where
\begin{eqnarray*}
(Q_{\varepsilon}u)(t)&=&(I-T(\omega))^{-1}\int_{t-\omega}^{t-\varepsilon}T_{\alpha}(t-s)F(s,u(s),u(s-\tau))ds\\
&=&T_{\alpha}(\varepsilon)(I-T(\omega))^{-1}\int_{t-\omega}^{t-\varepsilon}T_{\alpha}(t-s-\varepsilon)F(s,u(s),u(s-\tau))ds.
\end{eqnarray*}
 From Lemma 2.2, the operator $T_{\alpha}(\varepsilon)$ is compact in $X_{\alpha}$, it is follows that the set $(Q_{\varepsilon}\overline{\Omega}_{r})(t)$ is relatively compact in $X_{\alpha}$.  For any $u\in \overline{\Omega}_{r}$ and $t\in \R$, from the following inequality
\begin{eqnarray*}
&&\|Q_{1}u(t)-Q_{\varepsilon}u(t)\|_{\alpha}\\[8pt]
&\leq& C\int^{t}_{t-\varepsilon}\|T_{\alpha}(t-s)F(s,u(s),u(s-\tau))\|_{\alpha}ds\\ &\leq& C\int^{t}_{t-\varepsilon}\|A^{\alpha}T(t-s)F(s,u(s),u(s-\tau))\|ds\\
&\leq& CM_{\alpha}\int^{t}_{t-\varepsilon}\frac{h_{r}(s)}{(t-s)^{\alpha}}ds, \end{eqnarray*}
one can obtain that the set $(Q_{1}\overline{\Omega}_{r})(t)$ is relatively compact in $X_{\alpha}$ for all $t\in \R$.

Thus, the Arzela-Ascoli theorem guarantees that $Q_{1}$ is a compact operator.

Secondly, we prove that  $Q_{2}$ is a contraction. Let $u,v\in \overline{\Omega}_{r}$, by the condition (H2), Lemma 2.1(vi) and Lemma 2.6, we have
\begin{eqnarray*}
&&\|Q_{2}u(t)-Q_{2}v(t)\|_{\alpha}\\[8pt]
&=&\Big\|G(t,u(t-\xi))-(I-T(\omega))^{-1}\int^{t-\omega}_{t}T(t-s)AG(s,u(s-\xi))ds\\
&&-G(t,v(t-\xi))+(I-T(\omega))^{-1}\int^{t-\omega}_{t}T(t-s)AG(s,v(s-\xi))\Big\|_{\alpha}\\[8pt]
&\leq&\|G(t,u(t-\xi))-G(t,v(t-\xi))\|_{\alpha}\\[8pt]
&&+ \Big\|(I-T(\omega))^{-1}\int^{t-\omega}_{t}T(t-s)A(G(s,u(s-\xi))-G(s,v(s-\xi)))ds\Big\|_{\alpha}\\[8pt]
&\leq& \|A^{\alpha-1}(AG(t,u(t-\xi))-AG(t,v(t-\xi))\|\\[8pt]
&&+ \|(I-T(\omega))^{-1}\|\cdot\int^{t-\omega}_{t}\|AT(t-s)\|\cdot\|A(G(s,u(s-\xi))
-G(s,v(s-\xi)))\|ds\\[8pt]
&\leq&C_{1-\alpha}\|AG(t,u(t-\xi)-AG(t,v(t-\delta))\|\\
&&+ CM\cdot\int^{t-\omega}_{t}\|AG(s,u(s-\xi))-AG(s,v(s-\xi))\|ds\\
&\leq& C_{1-\alpha}L\|u(t-\xi)-v(t-\xi)\|_{\alpha}
+CM_{\alpha}L\int^{t-\omega}_{t}\frac{1}{(t-s)^{\alpha}}\|u(s-\xi)-v(s-\xi)\|_{\alpha}ds\\
&\leq&(C_{1-\alpha}L+CM_{\alpha}L\frac{\omega^{1-\alpha}}{1-\alpha})\|u-v\|_{C\alpha},
\end{eqnarray*}
therefore,
$$\|Q_{2}u-Q_{2}v\|_{C}\leq (C_{1-\alpha}L+CM_{\alpha}L\frac{\omega^{1-\alpha}}{1-\alpha})\|u-v\|_{C\alpha}.\eqno(3.9)$$
Since $CM_{\alpha}L\gamma+C_{1-\alpha}L+CM_{\alpha}L\frac{\omega^{1-\alpha}}{1-\alpha}<1$, so $C_{1-\alpha}L+CM_{\alpha}L\frac{\omega^{1-\alpha}}{1-\alpha}<1$, it follows that $Q_{2}$ is a contraction.

By Lemma 2.7, we know that $Q$ has a fixed point $u\in\overline{\Omega}_{r}$, that is, Eq (1.4) has a $\omega$-periodic mild solution. The proof is completed.$\square$
\vskip3mm
In the condition (H1), if the function $h_{r}$ is independent of $t$, we can easily obtain a constant $\gamma\geq0$ satisfying (H3). For example, we replace the condition (H1) with \vskip1mm
\noindent (H1$'$)  there are positive constants $a_{0},a_{1}$ and $K$ such that
$$\|F(t,x_{0},x_{1})\|\leq a_{0}\|x_{0}\|_{\alpha}+a_{1}\|x_{1}\|_{\alpha}+K$$
for $t\in\R$ and $x_{0},x_{1}\in X_{\alpha}$.

\vskip2mm
\noindent In this case, for any $r>0$ and $x_{0},x_{1}\in X_{\alpha}$ with $\|x_{0}\|_{\alpha},\|x_{1}\|_{\alpha}\leq r$, we have
$$\|F(t,x_{0},x_{1})\|\leq r(a_{0}+a_{1})+K:=h_{r}(t),\ \ \ t\in\R,$$
thus,
$$\liminf_{r\rightarrow\infty}\frac{1}{r}\int_{t-\omega}^{t}\frac{h_{r}(s)}{(t-s)^{\alpha}}ds=(a_{0}+a_{1})\frac{\omega^{1-\alpha}}{1-\alpha}:=\gamma>0.$$
Therefore, we have the following result
\vskip3mm
\noindent\textbf{Corollary 3.2} \emph{Let $A:D(A)\subset X\rightarrow X$ be a closed linear operator, and $-A$ generate a compact and exponentially stable  analytic operator semigroup $T(t)(t\geq0)$ in Banach space $X$. For $\alpha\in [0,1)$, we assume that  $G:\R\times X_{\alpha}\rightarrow X_{1}$ and $F : \R \times X_{\alpha}^{2} \rightarrow X$ are continuous
functions, and for every $x,x_{0},x_{1}\in X_{\alpha}$, $G(t,x),F(t,x_{0},x_{1})$ are $\omega$-periodic in $t$. If the  conditions (H1 $'$), (H2) and
\vskip2mm
\noindent (H3 $'$) $CM_{\alpha}(a_{0}+a_{1}+L)\frac{\omega^{1-\alpha}}{1-\alpha}+C_{1-\alpha}L<1$, where $C=\|(I-T(\omega))^{-1}\|$,
\vskip2mm
\noindent hold,
 then Eq. (1.4) has at least one $\omega$-periodic mild solution $u$.}

\vskip3mm

Furthermore, we assume that $F$ satisfies Lipschitz condition, namely,
\emph{\vskip1mm
\noindent (H1 $''$) there are positive constants $a_{0},a_{1}$, such that
$$\|F(t,x_{0},x_{1})-F(t,y_{0},y_{1})\|\leq a_{0}\|x_{0}-y_{0}\|_{\alpha}+a_{1}\|x_{1}-y_{1}\|_{\alpha}, \ \ \ t\in\R, x_{0},x_{1},y_{0},y_{1}\in X_{\alpha},$$}
then we can obtain the following result.
\vskip3mm
\noindent\textbf{Theorem 3.3.}\emph{ Let $A:D(A)\subset X\rightarrow X$ be a closed linear operator, and $-A$ generate a compact and exponentially stable  analytic operator semigroup $T(t)(t\geq0)$ in Banach space $X$. For $\alpha\in [0,1)$, we assume that  $G:\R\times X_{\alpha}\rightarrow X_{1}$ and $F : \R \times X_{\alpha}^{2} \rightarrow X$ are continuous
functions, and for every $x,x_{0},x_{1}\in X_{\alpha}$, $G(t,x),F(t,x_{0},x_{1})$ are $\omega$-periodic in $t$. If the conditions (H1 $''$),(H2), and (H3 $'$) hold, then Eq. (1.4) has unique $\omega$-periodic mild solution $u$.}
\vskip3mm
\noindent\textbf{Proof } From (H1$''$) we easily see that (H1$'$) holds. In fact, for any $t\in\R$ and $x_{0},x_{1}\in X_{\alpha}$, by the condition (H1$''$),
\begin{eqnarray*}\|F(t,x_{0},x_{1})\|&\leq&\|F(t,x_{0},x_{1})-F(t,\theta,\theta)\|+\|F(t,\theta,\theta)\|\\[8pt]
&\leq&a_{0}\|x_{0}\|_{\alpha}+a_{1}\|x_{1}\|_{\alpha}+\|F(t,\theta,\theta)\|.\end{eqnarray*}
From the continuity and periodicity of $F$, we can choose $K=\max_{t\in[0,\omega]}\|F(t,\theta,\theta)\|$, thus, the condition (H1$'$) holds. Hence by Corollary 3.2, Eq.(1.4) has $\omega$-periodic mild solutions. Let $u_{1},u_{2}\in C_{\omega}(\R,X_{\alpha})$ be the $\omega$ -periodic mild solutions of Eq.(1.4), then they are the fixed points of the operator $Q$ which is defined by (3.3). Hence,
\begin{eqnarray*}
&&\|Qu_{2}(t)-Qu_{1}(t)\|_{\alpha}\\[8pt]
&\leq&\Big\|(I-T(\omega))^{-1}\int_{t-\omega}^{t}T(t-s)\Big(F(s,u_{2}(s),u_{2}(s-\tau))-F(s,u_{1}(s),u_{1}(s-\tau))\Big)ds\Big\|_{\alpha}\\
&&+\Big\|G(t,u_{2}(t-\xi))-G(t,u_{1}(t-\xi))\Big\|_{\alpha}\\[8pt]
&&+\Big\|(I-T(\omega))^{-1}\int_{t-\omega}^{t}T(t-s)\Big(AG(s,u_{2}(s-\xi))-AG(s,u_{1}(s-\xi))\Big)ds\Big\|_{\alpha}\\
&\leq&  C\cdot\int_{t-\omega}^{t}\|A^{\alpha}T(t-s)\|\cdot \|F(s,u_{2}(s),u_{2}(s-\tau))-F(s,u_{1}(s),u_{1}(s-\tau)))\|ds\\[8pt]
&&+\|A^{\alpha-1}(AG(t,u_{2}(t-\xi))-AG(t,u_{1}(t-\xi)))\|\\[8pt]
&&+ C\cdot\int_{t-\omega}^{t}\|A^{\alpha}T(t-s)\|\cdot\|AG(s,u_{2}(s-\xi))-AG(s,u_{1}(s-\xi))\|ds\\ &\leq&CM_{\alpha}\cdot\int_{t-\omega}^{t}\frac{1}{(t-s)^{\alpha}}(a_{0}\|u_{2}(s)-u_{1}(s)\|_{\alpha}+a_{1}\|u_{2}(s-\tau)-u_{1}(s-\tau)\|_{\alpha})ds\\[8pt]
&&+C_{1-\alpha}L\|u_{2}(t-\xi)-u_{1}(t-\xi)\|_{\alpha}\\[8pt]
&&+ CM_{\alpha}\cdot\int_{t-\omega}^{t}\frac{L}{(t-s)^{\alpha}}\|u_{2}(t-\xi)-u_{1}(t-\xi)\|_{\alpha}ds\\
&\leq& CM_{\alpha}\frac{\omega^{1-\alpha}}{1-\alpha}(a_{0}+a_{1})\|u_{2}-u_{1}\|_{C\alpha}+C_{1-\alpha}L\|u_{2}-u_{1}\|_{C\alpha}
+CM_{\alpha}\frac{\omega^{1-\alpha}}{1-\alpha}L\|u_{2}-u_{1}\|_{C\alpha}\\[8pt]
&=&\Big(CM_{\alpha}(a_{0}+a_{1}+L)\frac{\omega^{1-\alpha}}{1-\alpha}+C_{1-\alpha}L\Big)\cdot\|u_{2}-u_{1}\|_{C_{\alpha}},
\end{eqnarray*}
which implies that $\|u_{2}-u_{1}\|_{C\alpha}=\|Qu_{2}-Qu_{1}\|_{C\alpha}\leq (CM\omega(a_{0}+a_{1}+L)+C_{1-\alpha}L)\cdot\|u_{2}-u_{1}\|_{C\alpha}$. From this and the condition (H3$'$), it follows that $u_{2}=u_{1}$. Thus, Eq.(1.4) has only one $\omega$-periodic mild solution.   $\square$

\section{Existence of Classical and Strong Solutions}

In this section, we discuss the regularity properties of the  $\omega$-periodic mild solution of Eq. (1.4), and present essential conditions on the nonlinearity $F$ and $G$ to guarantee that  Eq. (1.4) has $\omega$-periodic classical and strong solutions.

Now, we are in a position to state and prove the main result of this section:

\vskip3mm
\noindent\textbf{Theorem 4.1.} \emph{Let $A:D(A)\subset X\rightarrow X$ be a closed linear operator, and $-A$ generate a compact and exponentially stable  analytic operator semigroup $T(t)(t\geq0)$ in Banach space $X$. For $\alpha\in [0,1)$, we assume that  $G:\R\times X_{\alpha}\rightarrow X_{1}$ and $F : \R \times X_{\alpha}^{2} \rightarrow X$ are continuous
functions, and for every $x,x_{0},x_{1}\in X_{\alpha}$, $G(t,x),F(t,x_{0},x_{1})$ are $\omega$-periodic in $t$. If the following conditions
\vskip1mm
\noindent (H4) there exist $L_{1}$ and $\mu_{1}\in (0,1)$ such that $$\|F(t_{2},x_{0},,x_{1})-F(t_{1},y_{0},y_{1})\|\leq L_{1}(|t_{2}-t_{1}|^{\mu_{1}}+\|x_{0}-y_{0}\|_{\alpha}+\|x_{1}-y_{1}\|_{\alpha})$$
for each $t_{1},t_{2}\in \R$ and $x_{0},x_{1},y_{0},y_{1}\in X_{\alpha}$,
\vskip1mm
\noindent (H5) $G(t,\theta)=\theta$ for $t\in \R$,  there exist $L_{2}$ and $\mu_{2}\in (0,1)$ such that
$$\|AG(t_{2},x)-AG(t_{1},y)\|\leq L_{2}(|t_{2}-t_{1}|^{\mu_{2}}+\|x-y\|_{\alpha})$$
for each $t_{1},t_{2}\in \R$ and $x,y\in X_{\alpha}$,
\vskip1mm
\noindent (H6)  $CM_{\alpha}(2L_{1}+L_{2})\frac{\omega^{1-\alpha}}{1-\alpha}+C_{1-\alpha}L_{2}<1$, where $C=\|(I-T(\omega))^{-1}\|$,
\vskip1mm
\noindent hold, then Eq.(1.1) has an  $\omega$-periodic
classical solution.}
\vskip3mm
\noindent\textbf{Proof } Let $Q$ be the operator defined by (3.3) in the proof of
 Theorem 3.1. By the assumptions of Theorem 4.1 and the proof of Theorem 3.1, we know that the operator  $Q:C_{\omega}(\R,X_{\alpha})\rightarrow C_{\omega}(\R,X_{\alpha})$ is well defined. From the conditions (H4) and (H5), for any $u_{1},u_{2}\in C_{\omega}(\R,X_{\alpha}),t\in \R$, similar to the proof of theorem 3.3, we have
 \begin{eqnarray*}
\|Qu_{2}(t)-Qu_{1}(t)\|_{\alpha}\leq\Big(CM_{\alpha}(2L_{1}+L_{2})\frac{\omega^{1-\alpha}}{1-\alpha}+C_{1-\alpha}L_{2}\Big)\cdot\|u_{2}-u_{1}\|_{C_{\alpha}},
\end{eqnarray*}
which implies that
$$\|Qu_{2}-Qu_{1}\|_{C\alpha}\leq\Big(CM_{\alpha}(2L_{1}+L_{2})\frac{\omega^{1-\alpha}}{1-\alpha}+C_{1-\alpha}L_{2}\Big)\cdot\|u_{2}-u_{1}\|_{C_{\alpha}}<\|u_{2}-u_{1}\|_{C_{\alpha}}.\eqno(4.1)$$
Hence, $Q:C_{\omega}(\R,X_{\alpha})\rightarrow C_{\omega}(\R,X_{\alpha})$ is a contraction, thus $Q$  has unique fixed point $u_{0}\in C_{\omega}(\R,X_{\alpha})$. By the definition of $Q$, $u_{0}$ is  $\omega$-periodic mild solution of Eq. (1.4).

Next, we prove that $u_{0}$ is $\omega$-periodic classical solution. From the periodicity of $u_{0}$, we only need prove it on $[0,\omega]$. Let $t\in[0,\omega]$ and $h(t)=F(t,u(t),u(t-\tau))-AG(t,u(t-\delta))$, then $h\in C([0,\omega],X)$. For $\forall \varepsilon\in(0,\omega)$,  since $u_{0} $ is the $\omega$-periodic mild solution of Eq.(1.4), hence $u_{0} $ is the mild solution of the initial value problem
$$\left\{\begin{array}{ll}
(u(t)-G(t,u(t-\delta)))'+A(u(t)-G(t,u(t-\delta)))=h(t),\ t\in [\varepsilon,\omega],\\[8pt]
u(\varepsilon)=u_{0}(\varepsilon).
 \end{array} \right.\eqno(4.2)$$
While $u_{0}(\varepsilon)\in X_{\alpha}$, from Lemma 2.4, it  follows that
$$u_{0}\in C^{\mu_{3}}([\varepsilon,\omega],X_{\alpha-\mu_{3}})\hookrightarrow C^{\mu_{3}}([\varepsilon,\omega],X),\ \ \ \mu_{3}\in(0,\alpha) .$$
On the other hand, from the condition (H4) and (H5), we can deduce $h\in C^{\mu}([\varepsilon,\omega],X)$, where $\mu=\min\{\mu_{1},\mu_{2},\mu_{3}\}$. By Lemma 2.5, we obtain
that $u_{0}$ is classical solution of Eq. (4.2) and satisfies
$$u_{0}\in C^{1}((\varepsilon,\omega],X)\cap C([\varepsilon,\omega],X_{1}).$$
By the arbitrariness of $\varepsilon$, we claim that
$$u_{0}\in C^{1}([0,\omega],X)\cap C([0,\omega],X_{1}).$$
Therefore, $u_{0}$ is $\omega$-periodic classical solution of Eq.(1.4) and satisfies
$$u_{0}\in C_{\omega}^{1}(\R,X)\cap C_{\omega}(\R,X_{1}).$$
The proof is completed.$\square$

\vskip3mm
\noindent\textbf{Theorem 4.2.} \emph{Let $X$ be a reflexive Banach space, $A:D(A)\subset X\rightarrow X$ is a closed linear operator and $-A$ generates an exponentially stable and compact analytic semigroup $T(t)(t\geq0)$ in $X$.
 For $\alpha\in [0,1)$, we assume that  $G:\R\times X_{\alpha}\rightarrow X_{1}$ and $F : \R \times X_{\alpha}^{2} \rightarrow X$ are continuous
functions, and for every $x,x_{0},x_{1}\in X_{\alpha}$, $G(t,x),F(t,x_{0},x_{1})$ are $\omega$-periodic in $t$. If the conditions
\vskip1mm
\noindent (H4 $'$) there exists a constant $L_{1}>0$ such that
$$\|F(t_{2},x_{0},x_{1})-F(t_{1},y_{0},y_{1})\|\leq L_{1}(|t_{2}-t_{1}|+\|x_{0}-y_{0}\|_{\alpha}+\|x_{1}-y_{1}\|_{\alpha})$$
for any $t_{1},t_{2}\in \R$ and $x_{0},x_{1},y_{0},y_{1}\in X_{\alpha}$,
\vskip1mm
\noindent  (H5 $'$) $G(t,\theta)=\theta$ for $t\in \R$,  there exist $L_{2}$ and such that
$$\|AG(t_{2},x)-AG(t_{1},y)\|\leq L_{2}(|t_{2}-t_{1}|+\|x-y\|_{\alpha})$$
for each $t_{1},t_{2}\in \R$ and $x,y\in X_{\alpha}$,
 \vskip1mm
\noindent and (H6) hold, then Eq.(1.1) has an  $\omega$-periodic
strong solution $u$.}
\vskip3mm
\noindent\textbf{Proof } Let $Q$ be the operator defined by (3.3) in the proof of
 Theorem 3.1. For a given $r>0$, let $\overline{\Omega}_{r}\subset C_{\omega}(\R,X_{\alpha})$  is defined by (3.2). By the conditions (H4$'$-H6$'$), one can use the same argument as in the proof of Theorem 3.1 to obtain  that $(Q\overline{\Omega}_{r})\subset \overline{\Omega}_{r}$.

 For this $r$, consider the set
$$\overline{\Omega}=\{u\in C_{\omega}(\R,X_{\alpha})|\ \|u\|_{C\alpha}\leq r, \ \|u(t_{1})-u(t_{2})\|_{\alpha}<L^{*}|t_{2}-t_{1}|, t_{1},t_{2}\in \R\}                         \eqno(4.3)$$
for some $L^{*}$ large enough. It is clear that $\overline{\Omega}$ is convex closed and nonempty set. We shall prove that $Q$ has a fixed point on $\overline{\Omega}$. Obviously, from the proof of Theorem 3.1, it is sufficient to show that for any $u\in\overline{\Omega}$
$$\|(Qu)(t_{2})-(Qu)(t_{1})\|_{\alpha}\leq L^{*}|t_{2}-t_{1}|,\qquad \forall\ t_{1},t_{2}\in \R.\eqno(4.4)$$
In fact, by the definition of $Q$, the condition (H4 $'$),(H5 $'$) and (4.3), we have

\begin{eqnarray*}
&&\|Qu(t_{2})-Qu(t_{1})\|_{\alpha}\\[8pt]
&\leq&\Big\|(I-T(\omega))^{-1}\int_{t_{2}-\omega}^{t_{2}}T(t_{2}-s)F(s,u(s),u(s-\tau))ds\\
&&~~-(I-T(\omega))^{-1}\int_{t_{1}-\omega}^{t_{1}}T(t_{1}-s)F(s,u(s),u(s-\tau))\Big)ds\Big\|_{\alpha}\\
&&+\Big\|G(t_{2},u(t_{2}-\xi))-G(t_{1},u(t_{1}-\xi))\Big\|_{\alpha}\\
&&+\Big\|(I-T(\omega))^{-1}\int_{t_{2}-\omega}^{t_{2}}T(t_{2}-s)AG(s,u(s-\xi))ds\\[8pt]
&&~~-(I-T(\omega))^{-1}\int_{t_{1}-\omega}^{t_{1}}T(t_{1}-s)AG(s,u(s-\xi))ds\Big\|_{\alpha}\\[8pt]
&\leq&\Big\|(I-T(\omega))^{-1}\int_{0}^{\omega}A^{\alpha}T(s)\Big(F(t_{2}-s,u(t_{2}-s),u(t_{2}-s-\tau))\\[8pt]
&&-F(t_{1}-s,u(t_{1}-s),u(t_{1}-s-\tau))\Big)ds\Big\|\\[8pt]
&&+C_{1-\alpha}\|AG(t_{2},u(t_{2}-\xi))-AG(t_{1},u(t_{1}-\xi)\|\\[8pt]
&&+\Big\|(I-T(\omega))^{-1}\int_{0}^{\omega}A^{\alpha}T(s)\Big(AG(t_{2}-s,u(t_{2}-s-\xi))\\[8pt]
&&-AG(t_{1}-s,u(t_{2}-s-\xi))\Big)ds\Big\|\\[8pt]
&\leq&CM_{\alpha}\frac{\omega^{1-\alpha}}{1-\alpha}L_{1}(1+2L^{*})|t_{2}-t_{1}|
+C_{1-\alpha}L_{2}(1+L^{*})|t_{2}-t_{1}|\\[8pt]
&&+CM_{\alpha}\frac{\omega^{1-\alpha}}{1-\alpha}L_{2}(1+L^{*})|t_{2}-t_{1}|\\[8pt]
&=&\Big(CM_{\alpha}\frac{\omega^{1-\alpha}}{1-\alpha}(L_{1}+L_{2})
+C_{1-\alpha}L_{2}+(CM_{\alpha}(2L_{1}+L_{2})\frac{\omega^{1-\alpha}}{1-\alpha}
+C_{1-\alpha}L_{2})L^{*}\Big)|t_{2}-t_{1}|\\[8pt]
&:=&(K_{0}+K^{*}L^{*})|t_{2}-t_{1}|,\end{eqnarray*}
where $K_{0}=CM_{\alpha}\frac{\omega^{1-\alpha}}{1-\alpha}(L_{1}+L_{2})
+C_{1-\alpha}L_{2}$ is a constant independence of $L^{*}$, and $K^{*}=(CM_{\alpha}(2L_{1}+L_{2})\frac{\omega^{1-\alpha}}{1-\alpha}
+C_{1-\alpha}L_{2})<1$. Hence,
$$\|Qu(t_{1})-Qu(t_{1})\|_{\alpha}\leq L^{*}|t_{2}-t_{1}|, \ \ \ \mathrm{\ for\ all}\ t_{2},t_{1}\in \R,\eqno(4.5)$$
whenever $L^{*}\geq \frac{K_{0}}{1-K^{*}}$. Therefore, $Q$ has a fixed point $u$ which is an $\omega$-periodic mild solution of Eq.(1.4).

By the above calculation, we see that for this $u(\cdot)$, all the following functions
\begin{eqnarray*}&&g(t)=G(t,u(t-\xi)),\\[8pt]
&&\Phi(t)=(I-T(\omega))^{-1}\int_{t-\omega}^{t}T(t-s)F(s,u(s),u(s-\tau))ds,\\
&&\Psi(t)=(I-T(\omega))^{-1}\int_{t-\omega}^{t}T(t-s)AG(s,u(x-\xi))ds
\end{eqnarray*}
are Lipschitz continuous, respectively. Since the $u$ is Lipschitz continuous on $\R$ and the space $X_{\alpha}$ is reflexive by the assumption and Lemma 2.3, then a result of \cite{Komura1969} asserts that $u(\cdot)$ is a.e. differentiable on $\R$ and $u'(\cdot)\in L_{loc}^{1}(\R,X_{\alpha})$. Furthermore, by a stantdard arguement as  Theorem 4.2.4 in \cite{Pazy1993}, we can obtain that
\begin{eqnarray*}\Phi'(t)&=&(I-T(\omega))^{-1}\Big((I-T(\omega))F(t,u(t),u(t-\tau))\\
&&-\int_{t-\omega}^{t}AT(t-s)F(s,u(s),u(s-\tau))ds\Big)\\
\Psi'(t)&=&(I-T(\omega))^{-1}\Big((I-T(\omega))AG(t,u(t-\xi))\\
&&-\int_{t-\omega}^{t}AT(t-s)AG(s,u(x-\xi))ds\Big).
\end{eqnarray*}

Hence, for almost every $t\in \R$
\begin{eqnarray*}
u'(t)&=& \Phi'(t)+g'(t)-\Psi'(t)\\[8pt]
&=&(I-T(\omega))^{-1}\Big((I-T(\omega))F(t,u(t),u(t-\tau))\\
&&-\int_{t-\omega}^{t}AT(t-s)F(s,u(s),u(s-\tau))ds\Big)+G'(t,u(t-\xi))\\[8pt]
&&-(I-T(\omega))^{-1}\Big((I-T(\omega))AG(t,u(t-\xi))\\
&&-\int_{t-\omega}^{t}AT(t-s)AG(s,u(x-\xi))ds\Big)\\
&=&F(t,u(t),u(t-\tau))-AG(t,u(t-\xi))+G'(t,u(t-\xi))\\[8pt]
&&-A\Big((I-T(\omega))^{-1}\int_{t-\omega}^{t}T(t-s)F(s,u(s),u(s-\tau))ds\\
&&-(I-T(\omega))^{-1}\int_{t-\omega}^{t}T(t-s)AG(s,u(x-\xi))ds\Big)\\
&=&F(t,u(t),u(t-\tau))-AG(t,u(t-\xi))+G'(t,u(t-\xi))\\[8pt]
&&-A(u(t)-G(t,u(t-\xi))),
\end{eqnarray*}
which implies that
$$(u(t)-G(t,u(t-\xi)))'+Au(t)=F(t,u(t),u(t-\tau)), \qquad\mathrm{a.e}\quad  t\in\R.\eqno(4.6)$$
This shows that $u$ is  a strong solution for Eq.(1.4) and the proof is completed. $\square$\vskip3mm

\section{Application}
In this section, we present one example, which
indicates how our abstract results can be applied to concrete problems.

\vskip3mm
Consider the following time $\omega$-periodic solutions of the parabolic boundary value problem with delays
\begin{eqnarray*}  \left\{\begin{array}{ll}
\frac{\partial }{\partial t}(u(x,t)-g(x,t)(u(x,t-\xi)+\frac{\partial}{\partial x}u(x,t-\xi)))-\frac{\partial^{2}}{\partial x^{2}} u(x,t)\\[8pt]
=
 f(x,t,u(x,t),\frac{\partial}{\partial x}u(x,t),u(x,t-\tau),\frac{\partial}{\partial x}u(x,t-\tau)),\ \ \ \ x\in[0,1],\ t\in\R, \\[10pt]
u(0,t)=u(1,t)=0,
 \end{array} \right.\qquad\ (5.1)\end{eqnarray*}
where $g\in C^{2,1}([0,1]\times\R)$, $g(0,\cdot)=g(1,\cdot)=0$, and $f\in C([0,1]\times\R\times\R\times\R\times\R\times\R)$, moreover $f,g$ are $\omega$-periodic in  second value, $\xi,\tau$ are positive constants which denote the time delays.

To treat this system in the abstract form (1.4), we choose the space  $X=L^{2}([0,1],\R)$, equipped with the $L^{2}$-norm $\|\cdot\|_{L^{2}}$, thus, $X$ is reflexive.

Define  operator $A:D(A)\subset X\rightarrow X$ by $$D(A):=\{u\in X|\ u'',u'\in X, u(0)=u(1)=0\},\quad Au=-\frac{\partial^{2}u}{\partial x^{2}}.\eqno(5.2)$$
Then $-A$ generates an exponentially stable compact analytic semigroup $T(t) (t \geq 0)$ in $X$.
It is well known that  $0\in\rho(A)$ and so the fractional powers of $A$ are well defined.
Moveover, $A$ has a discrete spectrum with eigenvalues of the form $n^{2}\pi^{2},\ n\in \mathbb{N}$, and the associated normalized eigenfunctions are given by $e_{n}(x)=\sqrt{2}\sin(n\pi x)$ for $x\in[0,1]$, the associated semigroup $T(t)(t\geq0)$ is explicitly given by
$$T(t)u=\sum_{n=1}^{\infty}e^{-n^{2}\pi^{2}t}( u, e_{n}) e_{n},  \ \ t\geq0,u\in X, \eqno(5.3)$$
where $(\cdot,\cdot)$ is an inner product on $X$, and it is not difficult to
verify that $\|T(t)\|\leq e^{-\pi^{2}t}$ for all $t\geq 0$.
Hence, we take $M=1$, $M_{\frac{1}{2}}=\Gamma(\frac{1}{2})$ and $\|I-T(\omega)\|\leq \frac{1}{1-e^{-\pi^{2}\omega}}$. The following results are
 also well known.\\
(e1) If $u\in D(A)$ then
$$Au=\sum_{n=1}^{\infty}n^{2}\pi^{2}( u,e_{n}) e_{n}. $$
(e2) For each $u\in X$,
$$ A^{-\frac{1}{2}}u=\sum_{n=1}^{\infty}\frac{1}{n}( u,e_{n}) e_{n}.$$
(e3) For each  $u\in D( A^{\frac{1}{2}}):=\{u\in X| \ \sum_{n=1}^{\infty}n(u,e_{n}) e_{n}\in X\}$,
$$ A^{\frac{1}{2}}u=\sum_{n=1}^{\infty}n(u,e_{n})e_{n},$$
 and $\|A^{-\frac{1}{2}}\|=1$.
\vskip3mm
The proof of the following lemma can be found in \cite{Travis1987}.
\vskip3mm
\noindent\textbf{Lemma 5.1.} \emph{If $v\in D(A^{\frac{1}{2}})$, then $v$ is absolutely continuous with $v'\in X$ and $\|v'\|_{L^{2}}=\|A^{\frac{1}{2}}v\|_{L^{2}}$.}
\vskip3mm
According to Lemma 5.1, we define the Banach space $X_{\frac{1}{2}}:=(D(A^{\frac{1}{2}}),\|\cdot\|_{\frac{1}{2}})$, where $\|v\|_{\frac{1}{2}}=\|A^{\frac{1}{2}}v\|_{L^{2}}$ for all $v\in X_{\frac{1}{2}}$.
\vskip3mm

Define
\begin{eqnarray*}F(t,u(t),u(t-\tau))(x)
&=&f(x,t,u(x,t),\frac{\partial }{\partial x}u(x,t),u(x,t-\tau),\frac{\partial }{\partial x}u(x,t-\tau)),\\
G(t,u(t-\xi))(x)&=&g(x,t)(u(x,t-\xi)+\frac{\partial }{\partial x}u(x,t-\xi)).\end{eqnarray*}
It is clear that $G:\R\times X_{\frac{1}{2}}\rightarrow X_{1}$ and $F:\R\times X_{\frac{1}{2}}\times X_{\frac{1}{2}}\rightarrow X$.
Then the partial differential equation with delays (5.1)
can be rewritten into the abstract evolution equation with delays (1.4).

\vskip3mm
\noindent\textbf{Theorem 5.1.} If the following conditions
\vskip1mm
\noindent(F1) there are positive constants $a_{0},a_{1}$ and $K$ such that for every $(x,t,v,\eta,w,\zeta)\in[0,1]\times\R\times\R\times\R\times\R\times\R$
$$|f(x,t,v,\eta,w,\zeta)|\leq a_{0}(|v|+|\eta|)+a_{1}(|w|+|\zeta|)+K,$$
\vskip1mm
\noindent(F2) there exist constants $L>0$  such that for every $t,v_{i},\eta_{i}\in\R(i=1,2)$,
$$|\frac{\partial^{2}}{\partial x^{2}}g(x,t)(v_{2}+\eta_{2})-\frac{\partial^{2}}{\partial x^{2}}g(x,t)(v_{1}+\eta_{1})|\leq L(|v_{2}-v_{1}|+|\eta_{2}-\eta_{1}|),$$
\vskip1mm
\noindent (F3) $\frac{2\omega^{\frac{1}{2}}}{1-e^{-\pi^{2}\omega}}\Gamma(\frac{1}{2})(a_{0}+a_{1}+L)+L<\frac{\pi}{1+\pi}$,
\vskip1mm
\noindent
hold, then the neutral partial differential equation with delays (5.1) has at least one time $\omega$-periodic mild solution.
\vskip3mm
\noindent\textbf{Proof.} Let $\phi,\varphi\in X_{\frac{1}{2}}$, from the condition (F1), we can get
\begin{eqnarray*}
&&\|F(t,\phi,\varphi)\|_{L^{2}}\\[8pt]
&=&\Big(\int^{1}_{0}\Big(f(x,t,\phi(x,t),\frac{\partial }{\partial x}\phi(x,t),\varphi(x,t),\frac{\partial }{\partial x}\varphi(x,t))\Big)^{2}dx\Big)^{\frac{1}{2}}\\
&\leq& \Big(\int_{0}^{1}\Big(a_{0}|\phi(x,t)+\frac{\partial }{\partial x}\phi(x,t)|+a_{1}|\varphi(x,t)+\frac{\partial }{\partial x}\varphi(x,t)|+K\Big)^{2}dx\Big)^{\frac{1}{2}}\\[8pt]
&\leq& a_{0}(\|\phi\|_{L^{2}}+\|\phi'\|_{L^{2}})+a_{1}(\|\varphi\|_{L^{2}}
+\|\varphi'\|_{L^{2}})+K\\[8pt]
&\leq&a_{0}(\frac{1}{\pi}+1)\|\phi'\|_{L^{2}}
+a_{1}(\frac{1}{\pi}+1)\|\varphi'\|_{L^{2}}+K\\[8pt]
&=&a_{0}(\frac{1}{\pi}+1)\|\phi\|_{\frac{1}{2}}
+a_{1}(\frac{1}{\pi}+1)\|\varphi\|_{\frac{1}{2}}+K,
\end{eqnarray*}
thus, the condition (H1$'$) in Section 3 holds.

Let $\phi,\varphi\in X_{\frac{1}{2}}$, from the condition (F2),  we have
 \begin{eqnarray*}
 &&\|AG(t,\phi)-AG(t,\varphi)\|_{L^{2}}\\
 &=&\left(\int^{1}_{0}\left(\frac{\partial^{2}}{\partial x^{2}}g(x,t)\cdot\Big((\phi(x,t)+\frac{\partial}{\partial x}\phi(x,t))-(\varphi(x,t)+\frac{\partial}{\partial x}\varphi(x,t))\Big)\right)^{2}dx\right)^{\frac{1}{2}}\\
 &\leq& \Big(\int^{1}_{0}L^{2}\Big(|\phi(x,t)-\varphi(x,t)|+|\frac{\partial}{\partial x}\phi(x,t)-\frac{\partial}{\partial x}\varphi(x,t)|\Big)^{2}dx\Big)^{\frac{1}{2}}\\[8pt]
 &\leq& L(\|\phi-\varphi\|_{L^{2}}+\|(\phi-\varphi)'\|_{L^{2}})\\[8pt]
 &\leq&L(\frac{1}{\pi}+1)\|\phi-\varphi\|_{\frac{1}{2}},
 \end{eqnarray*}
thus, the condition (H2) in  Section 3 holds.

Finally, by (F3), we can easily to prove that the condition (H3 $'$) holds in Section 3.
Therefore, from Corollary 3.2, it follows that the neutral partial differential equation with delays (5.1) has at least one time $\omega$-periodic mild solution.
The proof is completed. $\square$
\vskip3mm

For showing the existence of classical and strong solutions, the following assumptions are need:
\vskip3mm
\noindent
(F4)  there exist constants $l_{1}$ and $\mu_{1}\in(0,1]$ such that for every $t_{i},v_{i},w_{i},\eta_{i},\zeta_{i}\in\R$,
\begin{eqnarray*}&&|f(x,t_{2},v_{2},\eta_{2},w_{2},\zeta_{2})-f(x,t_{1},v_{1},\eta_{1},w_{1},\zeta_{1})|\\
&\leq& l_{1}(|t_{2}-t_{1}|^{\mu_{1}}+|v_{2}-v_{1}|+|\eta_{2}-\eta_{1}|+|w_{2}-w_{1}|+|\zeta_{2}-\zeta_{1}|),\ \ \ x\in[0,1],\end{eqnarray*}
 \vskip1mm
\noindent
(F5)  there exist constants $l_{2}>0$ and $\mu_{2}\in(0,1]$ such that for every $t_{i},v_{i},\eta_{i}\in\R(i=1,2)$,
$$|\frac{\partial^{2}}{\partial x^{2}}g(x,t_{2})(v_{2}+\eta_{2})-\frac{\partial^{2}}{\partial x^{2}}g(x,t_{1})(v_{1}+\eta_{1})|\leq l_{2}(|t_{2}-t_{1}|^{\mu_{2}}+|v_{2}-v_{1}|+|\eta_{2}-\eta_{1}|),,\ \ \ x\in[0,1],$$
\vskip1mm
\noindent(F6) $\frac{2\omega^{\frac{1}{2}}}{1-e^{-\pi^{2}\omega}}\Gamma(\frac{1}{2})(2l_{1}+l_{2})+l_{2}<\frac{\pi}{1+\pi}$.
\vskip3mm
\noindent Hence, for every $t_{1},t_{2}\in\R$ and $\phi_{1},\varphi_{1},\phi_{2},\varphi_{2}\in X_{\frac{1}{2}}$, we have
\begin{eqnarray*}
&&\|F(t_{2},\phi_{2},\varphi_{2})-F(t_{1},\phi_{1},\varphi_{1})\|_{L^{2}}\\[8pt]
&=&\Big(\int^{1}_{0}\Big(f(x,t_{2},\phi_{2}(x,t),\frac{\partial }{\partial x}\phi_{2}(x,t),\varphi_{2}(x,t),\frac{\partial }{\partial x}\varphi_{2}(x,t))\\
&&-f(x,t_{1},\phi_{1}(x,t),\frac{\partial }{\partial x}\phi_{1}(x,t),\varphi_{1}(x,t),\frac{\partial }{\partial x}\varphi_{1}(x,t))\Big)^{2}dx\Big)^{\frac{1}{2}}\\[8pt]
&\leq& \Big(\int_{0}^{1}l^{2}\Big(|t_{2}-t_{1}|^{\mu_{1}}+|\phi_{2}(x,t)-\phi_{1}(x,t)|+|\frac{\partial }{\partial x}\phi_{2}(x,t)-\frac{\partial }{\partial x}\phi_{1}(x,t)|\\[8pt]
&&+|\varphi_{2}(x,t)-\varphi_{1}(x,t)|+|\frac{\partial }{\partial x}\varphi_{2}(x,t)-\frac{\partial }{\partial x}\varphi_{1}(x,t)|\Big)^{2}dx\Big)^{\frac{1}{2}}\\[8pt]
&\leq& l_{1}(|t_{2}-t_{1}|^{\mu_{1}}+\|\phi_{2}-\phi_{1}\|_{L^{2}}
+\|(\phi_{2}-\phi_{1})'\|_{L^{2}}+\|\varphi_{2}-\varphi_{1}\|_{L^{2}}
+\|(\varphi_{2}-\varphi_{1})'\|_{L^{2}})\\[8pt]
&\leq& l_{1}(|t_{2}-t_{1}|^{\mu_{1}}+(\frac{1}{\pi}+1)\|(\phi_{2}-\phi_{1})'\|_{L^{2}}
+(\frac{1}{\pi}+1)\|(\varphi_{2}-\varphi_{1})'\|_{L^{2}})\\[8pt]
&\leq& l_{1}(\frac{1}{\pi}+1)(|t_{2}-t_{1}|^{\mu_{1}}+
\|(\phi_{2}-\phi_{1})\|_{\frac{1}{2}}+\|\varphi_{2}-\varphi_{1}\|_{\frac{1}{2}}),
\end{eqnarray*}
and
\begin{eqnarray*}
 &&\|AG(t_{2},\phi_{2})-AG(t_{1},\phi_{1})\|_{L^{2}}\\[8pt]
 &=&\Big(\int^{1}_{0}\Big(\frac{\partial^{2}}{\partial x^{2}}g(x,t_{2})(\phi_{2}(x,t_{2})+\frac{\partial}{\partial x}\phi_{2}(x,t_{2}))\\
 &&-\frac{\partial^{2}}{\partial x^{2}}g(x,t_{1})(\phi_{1}(x,t_{1})+\frac{\partial}{\partial x}\phi_{1}(x,t_{1}))\Big)^{2}dx\Big)^{\frac{1}{2}}\\
 &\leq& \Big(\int^{1}_{0}l_{2}^{2}\Big(|t_{2}-t_{1}|^{\mu_{2}}
 +|\phi_{2}(x,t_{2})-\phi_{1}(x,t_{1})|+|\frac{\partial}{\partial x}\phi_{2}(x,t_{2})-\frac{\partial}{\partial x}\phi_{1}(x,t_{1})|\Big)^{2}dx\Big)^{\frac{1}{2}}\\[8pt]
 &\leq& l_{2}(|t_{2}-t_{1}|^{\mu_{2}}+\|\phi_{2}-\phi_{1}\|_{L^{2}}
 +\|(\phi_{2}-\phi_{1})'\|_{L^{2}})\\[8pt]
 &\leq&l_{2}(\frac{1}{\pi}+1)(|t_{2}-t_{1}|^{\mu_{2}}
 +\|\phi-\varphi\|_{\frac{1}{2}}),
 \end{eqnarray*}
which implies that the conditions (H4) and (H5) for $\mu_{1},\mu_{2}\in(0,1)$ (or (H4$'$) and (H5$'$) for $\mu_{1}=\mu_{2}=1$) hold. On the other hand, by the condition (F6), we can easily prove that the condition (H6) holds in Section 4.

Consequently, all the conditions stated in Theorem 4.1 and Theorem 4.2 are satisfied and we obtain the following interesting results.

\vskip3mm
\noindent\textbf{Theorem 5.2.} If the conditions (F4-F5) hold for $\mu_{1},\mu_{2}\in(0,1)$, then the neutral partial differential equation with delays (5.1) exists time $\omega$-periodic classical solution.
\vskip3mm
\noindent\textbf{Theorem 5.3.} If the conditions (F4-F5) hold for $\mu_{1}=\mu_{2}=1$, then the neutral partial differential equation with delays (5.1) exists time $\omega$-periodic strong solution.

\bibliographystyle{abbrv}

\begin{thebibliography}{99}
\small
\setlength{\parskip}{0pt}
\setlength{\baselineskip}{8pt}
\vspace{1pt}




\bibitem{Hale1993}{J.K. Hale, S.M. Verduyn Lunel, Introduction to Functional Differential Equations, Springer-Verlag, Berlin, 1993.}
\bibitem{Wu1996}{J. Wu, Theory and Application of Partial Functional Differential Equations, Springer-Verlag, New York, 1996.}
\bibitem{Wu1996}{J. Wu, H. Xia, Self-sustained oscillations in a ring array of coupled lossless transmission lines, J. Differ. Equ. 124 (1996) 247-278.}
\bibitem{Wu1999}{J. Wu, H. Xia, Rotating waves in neutral partial functional differential equations, Journal of Dynamics and Differential Equations, 11 (1999) 209-238.}
\bibitem{Adimy2004}{M. Adimy, H. Bouzahir, K. Ezzinbi, Existence and stability for some partial neutral functional differential equations with infinite delay, J. Math. Anal. Appl., 294 (2004) 438-461.}
\bibitem{Adimy2006}{M. Adimy, K. Ezzinbi, Existence and stability in the ¦Á-norm for partial functional equations of neutral type, Ann. Mat. Pura Appl., 185 (2006) 437-460.}
\bibitem{Ezzinbi2004}{ K. Ezzinbi, X.L. Fu, Existence and regularity of solutions for some neutral partial differential equations with non-local conditions, Nonlinear Anal., 57 (2004) 1029-1041}
\bibitem{Ezzinbi2010}{K. Ezzinbi, S. Ghnimib, Existence and regularity of solutions for neutral partial functional integrodifferential equa-tions, Nonlinear Anal. Real World Appl., 11 (2010) 2335-2344.}
\bibitem{Hernandeza2011}{ E. Hernandeza, M. Pierri, A. Prokopczyk, On a class of abstract neutral functional differential equations, Nonlinear Anal., 74 (2011) 3633-3643.}
\bibitem{Regan2011}{E. Hernandez, D. O'Regan, On a new class of abstract neutral differential equations, J. Funct. Anal., 261 (2011) 3457-3481.}
\bibitem{Adimy1998}{   M. Adimy and K. Ezzinbi, A class of linear partial neutral functional differential equations
with nondense domain, J. Differ. Equ., 147(1998)285-332.}

\bibitem{Babram1996}{ M.A. Babram, K. Ezzinbi, Periodic solutions of functional differential equations of neutral type, J. Math. Anal. Appl., 204 (1996) 898-909.}
\bibitem{Hernandez1998}{ E. Hernandez, H.R. Henriquez, Existence of periodic solutions of partial neutral functional differential equations with unbounded delay, J. Math. Anal. Appl., 221 (1998) 499-522.}


\bibitem{Benkhalti2006}{R. Benkhalti,  A. Elazzouzi, K. Ezzinbi, Periodic solutions for some partial neutral functional differential equations. Electron. J. Differ. Equ., 2006(2006) No.56, pp. 1-14.}
\bibitem{Fu2007}{ X. Fu, X. Liu, Existence of periodic solutions for abstract neutral non-autonomous equations with infinite delay. J. Math. Anal. Appl., 325 (2007) 249-267.}

\bibitem{Fu2008}{X. Fu, Existence of solutions and periodic solutions for abstract neutral equations with unbounded delay. Dyn. Contin. Discrete Impuls. Syst. Ser. A Math. Anal., 15 (2008) 17-35.}
\bibitem{Benkhalti2010}{ R. Benkhalti, A. Elazzouzi, K. Ezzinbi, Periodic solutions for some nonlinear partial neutral functional differential equations. Internat. J. Bifur. Chaos Appl. Sci. Engrg., 20 (2010) 545-555.}
\bibitem{Ezzinbi2014}{K. Ezzinbi, B. A. Kyelem, S. Ouaro, Periodicity in the $\alpha$-norm for partial functional differential equations in fading memory spaces. Nonlinear Anal., 97 (2014) 30-54.}




\bibitem{Burton1991}{T.A. Burton, B. Zhang, Periodic solutions of abstract differential equations with infinite delay, J. Diffe. Equ.,
90 (1991) 357-396.}

\bibitem{Xiang1992}{X. Xiang, N.U. Ahmed, Existence of periodic solutions of semilinear evolution equations with time lags, Nonlinear
Anal., 18 (1992) 1063-1070.}

\bibitem{Liu1998}{J. Liu, Bounded and periodic solutions of finite delays evolution equations, Nonlinear Anal., 34 (1998) 101-111.}

\bibitem{Liu2000}{J. Liu, Periodic solutions of infinite delay evolution equations, J. Math. Anal. Appl., 247 (2000) 644-727.}

\bibitem{Liu2003}{J. Liu, Bounded and periodic solutions of infinite delay evolution equations, J. Math. Anal. Appl., 286 (2003) 705-
712.}
\bibitem{Massera1950}{J. Massera, The existence of periodic solutions of differential equations, Duke Math. J., 17(1950)426-429. }

 \bibitem{Zhu2008}{J.  Zhu, Y. Liu, Z. Li, The existence and attractivity of time periodic solutions for evolution equations with
delays, Nonlinear Anal. (Real World Appl.), 9 (2008) 842-851.}
\bibitem{Li2011}{Y. Li, Existence and asymptotic stability of periodic solution for evolution equations with delays, J. Functional Anal., 261 (2011) 1309-1324.}
\bibitem{Pazy1993}{A. Pazy, Semigroup of linear operators and applications to partial differential equations, Springer-Verlag, New York, 1993.}
\bibitem{Triggiani1975}{R. Triggiani, On the stabilizability problem in Banach space, J. Math. Anal. Appl., 52 (1975) 383-403.}
 \bibitem{Hino1987}{Y. Hino, S. Murakami, Periodic solutions of linear voltera systems in differential equations, in: Lecture Notes in Pure and App. Math., Vol. 118, Dekker, VY, 1987, pp. 319-326.}

\bibitem{Liu2009}{H. Liu, J. Chang, Existence for a class of partial differential equations with nonlocal conditions, Nonlinear Anal., 70(2009) 3076-3083.}
\bibitem{Chang2009}{J. Chang, H. Liu, Existence of solutions for a class of neutral partial differential equations with nonlocal conditions in the $\alpha$-norm,  Nonlinear Anal., 71(2009) 3759-3768.}

\bibitem{Li2005}{Y. Li, Existence and uniquness of positive periodic solution for abstract semilinear evolution equations, J.Sys. Sci. Math. Scis., 25 (2005) 720-728 (in Chinese). }
\bibitem{Sadovskii1967}{B.N. Sadovskii, On a fixed point principle, Funct. Anal. Appl.,  1 (1967) 74-76.}

\bibitem{Komura1969}{J. Komura, Differentiability of nonlinear semigroups, J. Math. Soc. Japan, 21(1969)375-402.}
\bibitem{Travis1987}{C.C. Travis, G.F. Webb, Existence, stability and compactness with $\alpha$-norm for partial functional differential equations, Transl. Amer. Math. Soc., 240
(1978) 129-143.}




\end{thebibliography}

\end{document}